\documentclass{article}
\usepackage{xypic}
\usepackage{amsthm}
\LaTeXdiagrams
\usepackage{graphicx}

\usepackage{amssymb}

\newcommand{\bP}{\Bbb P}
\newcommand{\C}{\mathbb C}
\newcommand{\Q}{\Bbb Q}
\newcommand{\Z}{\Bbb Z}
\newcommand{\R}{\Bbb R}

\newtheorem{theorem}{Theorem}
\newtheorem{definition}{Definition}
\newtheorem{proposition}{Proposition}
\newtheorem{lemma}{Lemma}
\newtheorem{corollary}{Corollary}

\theoremstyle{remark}

\begin{document}

\title{The Chow ring of the Classifying Space\\$BSO(2n,\C)$}
\author{Rebecca Field}
\date{March 22, 2001 \footnote{I'd like to apologize for not posting this sooner.}}
\maketitle

{\footnotesize  
\begin{quote}
ABSTRACT. We compute the Chow ring of the classifying space $BSO(2n,\C)$ 
in the sense of Totaro 
using the fibration $Gl(2n)/SO(2n) \to BSO(2n) \to BGl(2n)$ and a computation 
of the Chow ring of $Gl(2n)/SO(2n)$ in a previous paper.  We find this Chow 
ring is generated by Chern classes and a characteristic class defined by 
Edidin and Graham which maps to $2^{n-1}$ times the Euler class under the 
usual class map from the Chow ring to ordinary cohomology.  Moreover, we show this class represents $1/2^{n-1}(n-1)!$ times 
the $n^{th}$ Chern class of the representation of $SO(2n)$ whose highest 
weight 
vector is twice that of the half-spin representation.
\end{quote}
}

Throughout this paper we write $Gl(n)$, $O(n)$, $SO(n)$, etc. 
for the
the complex algebraic groups of these types.  They are homotopy 
equivalent to the compact groups $U(n)$, $O(n, \R )$, and $SO(n, \R )$.

This paper is devoted to computing the Chow ring of the classifying space $BSO(2n)$.  The definition of the Chow ring of a classifying space we use is that of Totaro \cite{Burt2}.  This is defined as the limit of Chow rings of finite dimensional algebraic varieties 
approximating the classifying space, and coincides with the ring of characteristic classes for principal 
$SO(2n,\C)$ bundles over smooth varieties.

Our main result is:

\begin{theorem}
$$CH^*BSO(2n) \cong \Z[c_2, ..., c_{2n}, y_n] / 
( 2c_{odd}, y_n.c_{odd}, y^2_n + (-1)^n 2^{2n-2}c_{2n} ),$$
where $y_n$ 
is the class defined by Edidin and Graham \cite{E&G} which 
maps to $2^{n-1}\chi$ in cohomology
under the usual class map.
\end{theorem}

For $n\leq 3$,
this Chow ring is
isomorphic to the quotient of complex cobordism: 
$MU^*BSO(2n)/MU^{<0}(pt).MU^*BSO(2n)$ via
Totaro's cycle map \cite{Burt1}.  
(This is also the case for $Gl(n)$, $SL(n)$, $O(n)$, 
and $SO(2n+1)$.)  Totaro has conjectured that if the 
cohomology of the classifying space is concentrated in even degrees then 
the Chow ring is isomorphic to this quotient of complex cobordism.
For $n>3$, however, the complex cobordism of
$BSO(2n)$ is not known.  

Theorem 1 is proven via the fibration $Gl(2n)/SO(2n) \to BSO(2n) \to BGl(2n)$, and so depends on the following computation from  
{\em The Chow ring of the symmetric space $Gl(2n)/SO(2n)$}\cite{me}.

\begin{theorem}
$$CH^*(Gl(2n)/SO(2n)) \cong \Z x_0 \oplus \Z y$$ where $x_0$ is the 
codimension 
zero cycle and $y$ is in codimension $n$.
\end{theorem}

The Chow ring of $BSO(2n)$ is not generated by Chern classes
of the standard representation (the generator 
$y_n$ doesn't come from representations of $Gl(2n)$).  
This raises the he question: Is the Chow ring generated by Chern 
classes of any representations of $SO(2n)$?  (As is the case for $SO(2n+1)$ and 
$O(2n)$ \cite{Burt2} \cite{Rahul}.)  

$SO(2n)$ has 
pair of 
representations 
for which $SO(2n+1)$ and 
$O(2n)$ have no analog,
namely the two representations 
whose highest 
weight vectors are twice that of the half spin representations of 
$\mathfrak{so}_{2n}$.  
However, these representations do not suffice to generate the Chow ring.
For $n \geq 3$, only a finite index subgroup of the Chow ring is generated by Chern classes 
of representations. 

Specifically, 

\begin{theorem} Let $D^+_n$ be the irreducible 
representation of $SO(2n)$ whose 
highest weight vector is twice that of one of the two half spin 
representations of $\mathfrak{so}_{2n}$.  Then,
$c_n(D^+_n) = 2^{n-1}(n-1)!\chi + p(c_1,c_2,...,c_{2n})$, in cohomology, where 
\linebreak $p(c_1,c_2,...c_{2n})$ is a polynomial in the Chern classes of the 
standard representation of $SO(2n)$ and $\chi$ is the Euler class.
\end{theorem}

Namely, even in ordinary cohomology, the $n$-th Chern class of $D_n^+$ only generates 
a multiple of the 
Euler class and not the Euler class itself.  
In fact, this multiple is  
$\dsize\frac{1}{n}|\frak{W}|$ where $|\frak{W}|$ is the order of the 
Weyl group. 
In the Chow ring, where we do not have an actual Euler class, rather $y_n$ 
mapping to $2^{n-1}\chi$, 
the theorem becomes
$c_n(D_n^+)=(n-1)!y_n + p(c_1, \ldots, c_n)$.

As an additional corollary, we will produce a useful 
stratification of the space of vector bundles with structure group $SO(2n)$.  
As Edidin and Graham observe \cite{E&G}, a Zariski locally 
trivial vector bundle has a characteristic class in the Chow ring which 
maps to the classical Euler class in cohomology via the usual class map.  
That is, a Z.L.T. bundle 
has an actual Euler class in the Chow ring.  However, as observed in 
Theorem 1, a general vector bundle whose structure group reduces to $SO(2n)$ 
can only be expected to have a class in the Chow ring 
representing $2^{n-1}$ times the 
Euler class.  
The stratification produced in section 5 has the property that a 
bundle $V \to X$ 
is in the $k^{th}$ strata 
if and only if there is a class in $CH^*X$ which maps to $2^{k-1}\chi$.  
This stratification will be constructed geometrically. 

\medskip

There is much work in the literature on generalized cohomology of classifying spaces, see for example, \cite{examples}.  Chow theory fails to be a generalized cohomology theory, having no long exact sequence, however, the results of this paper fit well in that general framework.  There is also a family of results specifically computing $CH^*(BG)$ for various $G$ \cite{Burt2} \cite{Rahul} \cite{Vezzosi}.  

Specific algebraic geometric applications do exists for these results.  For example, Chow rings of classifying spaces are a basic building block in the
study of Godeaux-Serre varieties.  A Godeaux-Serre variety is the 
quotient of a complete
intersection by a free finite group action, and these have historically
been a rich source of examples in algebraic geometry.  Atiyah and
Hirzebruch used them to find torsion elements in cohomology which are
not in the Chow ring (counterexamples to an integral form of the Hodge
conjecture).  Godeaux-Serre varieties are also the examples that 
Totaro used to find torsion
elements of the Chow ring which map to zero in cohomology 
\cite{Burt1}.

$BSO(2n)$, 
classifies complex vector bundles with
quadratic form and trivial determinant line bundle.  Hence, this work 
is of interest when one
studies invariants of such quadratic bundles.  In addition, $CH^*BSO(2n)$ is 
the coefficient ring for $SO(2n)$-equivariant intersection theory as 
defined in \cite{equiv.int.th}.

Most intriguingly, there are various other constructions of classifying spaces as simplicial schemes in appropriate homotopy categories \cite{morel&voevodsky}.  It would be interesting to compare these more sophisticated constructions to the results of this paper.

\noindent{\em Acknowledgments}\\
I would like to thank Burt Totaro under whose direction the majority of the work in this paper was done as a thesis project in 2000.  In addition, I would like to thank J. Peter May and Madhav Nori for their assistance as surrogate thesis advisors.  I'd also like to thank Frank Sottile, and Bill Graham for helpful conversations and Sandy Kutin for much needed help with combinatorics.   Most especially, I'd like to thank Ian Grojnowski who helped to feed and water both me and this paper in many, many ways.

\section{Outline of the Proof}

The traditional method for calculating cohomological invariants of the
classifying space is to decompose it into fibrations of Grassmannians
and then use Schubert calculus. This is the method used in the examples of
\cite{Burt1} (Totaro),
where it works reasonably well, and the method used by
Pandharipande to compute the Chow 
ring of $SO(4)$\cite{Rahul}. Pandharipande's computation, though 
beautiful, is long
and difficult, and it was unclear whether it would generalize. 
The results of this paper show that, in fact, this technique cannot work 
without fundamental changes 
for $SO(2n)$, $n>2$, since that Chow ring is not 
generated by Chern classes of representations and hence not easily 
detectable using Schubert calculus.
Thus a new technique was needed. 

In section 2, we review basic notation and results necessary for this paper, including 
the following result fundamental to our chosen method of computation.  Namely, (Totaro)\cite{Burt2}
$CH^*BSO(2n)$ is generated 
as a module over $CH^*BGL(2n)$ by any set of elements of $CH^*BSO(2n)$ which 
restrict to generators of $CH^*(GL(2n)/SO(2n))$ as an Abelian group.    

$Gl(2n)/SO(2n)$ possesses such a generator, $y_n$ in codimension $n$, and in section 3 we show that it is the same as an element previously found by Edidin and Graham.  
The element $y_n$ maps to $2^{n-1}\chi$ in cohomology, where $\chi$ is 
the Euler class.  Such an element was constructed 
for even rank vector bundles
$V\to X$ with non-degenerate quadratic form whose structure group 
reduces to $SO(2n)$ \cite{E&G}.
In our case, the appropriate 
bundle 
over $Gl(2n)/SO(2n)$ on which to compute this class 
is the trivial vector bundle with the
quadratic form at a point $(g,\pm \sqrt{|g|}) \in GL(2n)/SO(2n)$ 
given by $<u,v>=u^t\cdot g\cdot v$.  

Altogether, this gives us that the Chow ring of 
$BSO(2n)$ is generated by
the Chern classes of the standard representation and Edidin and Graham's 
characteristic class.

Computation of the various
relations in section 4 gives a final answer of Theorem 1. 

Section 5, which is not used in the rest of the paper, contains the aforementioned stratification of the space of vector bundles with structure group $SO(2n)$

In section 6, we ask whether $y_n$ is a polynomial in Chern classes of 
representations of
$SO(2n)$.  The representation ring of $SO(2n)$ is
generated by exterior powers of the standard representation and one
additional irreducible representation $D_n^+$.  $D_n^+$ is the $+1$ 
eigenspace of the
Hodge star operator on $\Lambda^n(V)$, also known as the space of self-dual
$n$-forms.  Since these generate the representation ring, if $y_n$ is to be a
polynomial in the Chern classes of representations it must be
equal to $\pm c_n(D_n^+)$ mod Chern classes of the standard
representation.

By restricting to the maximal torus, we reduce the problem to 
that of finding a specific coefficient in a product of binomials.  
By a combinatorial argument due to 
S. Kutin we get Theorem 3 of the introduction, namely
the $n$-th 
Chern class of the representation $D_n^+$ only generates $(n-1)!y_n$ and 
not $y_n$ itself.   

These numbers suggested (as Bill Graham and Frank 
Sottile noticed) the action of the Weyl group 
enters into this computation 
(the order of the Weyl group is $2^{n-1}n!$, $8\cdot3$, $48\cdot4$,
and $384\cdot5$ for $n=3$, $4$ and $5$ respectively).  
To this end, Frank Sottile wrote a 
computer program computing the $d_n$ for $n$ up to $10$ and the 
order 
of the Weyl group continues to appear.  
This was later proven using Sandy Kutin's combinatorial argument.

\section{Basic notation and results}

For notation and conventions on Chow rings, we will refer to Fulton's 
{\em Intersection theory} \cite{Fulton}. 
The difference is only notational.  Fulton uses $A_*(X)$ to denote the ring of 
cycles mod rational equivalence, while we use the notation $CH_*(X)$.

The basic result about Chow rings that we will need are the following.

\begin{lemma} Let $X$ be a variety with 
a closed subvariety $Y$.  Then
$$CH_*Y \to CH_*X \to CH_*(X-Y) \to 0$$
is exact.
\end{lemma}

\begin{proposition} Let $X$ be a complex variety.  Then the class map 
$$\sc{cl}:CH^i(X)\to H^{2i}(X)$$
is a ring homomorphism.
\end{proposition}

Note that for $X$ a non-singular variety of dimension $n$, $CH_iX\cong CH^{n-i}X$.

The following theorem is the foundation for this method of computing the Chow 
ring of a classifying space.

\begin{theorem}(Totaro)\cite{Burt2}
Let $G$ be an algebraic group over a field, and let 
$G \hookrightarrow H$ be an embedding of $G$ into a group $H$ which 
is a product of the groups $GL(n_i)$ for some integers $n_i$.  Then the 
Chow ring $CH^*BG$ is generated as a module over $CH^*BH$ by elements 
of degree at most $dim (H/G)$.  This follows from the 
more precise statement:
$$CH^*(H/G) = CH^*BG \otimes_{CH^*BH} \Z.$$
\end{theorem}
Note $CH^*BH$ is a polynomial 
ring over $\Z$ generated by the Chern classes.

This theorem relies on the fact that the fibration

$$
\xymatrix{
H/G \ar[r]&BH \ar[d]\\
& BG\\
}
$$
\noindent is Zariski locally trivial.  
The result then follows from the fundamental exact sequence of 
Chow rings. 

We will apply this result when $G=SO(2n)$ and $H=Gl(2n)$.  To do this 
we need to know $CH^*(Gl(2n)/SO(2n))$ which is the main result of \cite{me}.

\begin{theorem} \cite{me} 
$CH^*(Gl(2n)/SO(2n)) \cong \Z \oplus \Z y_n$, where $y_n$ is a 
codimension $n$ cycle.
\end{theorem}

In the next section, we will review this computation to the extent necessary for 
comparison of this codimension $n$ element 
and Edidin and 
Graham's codimension $n$ characteristic class.

\section{The Euler class}

Define a {\em quadratic bundle} over a variety $X$ to be a pair $(V,q)$ where 
$V \to X$ is a vector bundle and $q$ is a non-degenerate quadratic form on V 
whose determinant is trivial.  

Let $X=Gl(2n)/SO(2n)$.  Then 
$X\cong\{(g,\epsilon )\mid g \in Gl(2n)\mbox{ symmetric, and }\\ 
\epsilon^2=det (g)\}$.  This is a double cover of $Gl(2n)/O(2n)$ 
with projection $\pi$.  
Let $V=X \times \C^{2n}$ be the trivial 
$2n$-dimensional vector bundle over $X$.  This admits a quadratic form 
whose value on $V_{(g,\epsilon)}$ is given by $<v,w>=v^tgw$.  This is a
universal bundle in the following sense: 
if $\mathcal{V} \to Y$ is a quadratic bundle such that 
$\mathcal{V}$ is trivial as a vector bundle, i.e. 
$\mathcal{V} = Y \times \C^{2n}$, then there is a map $f:Y \to X$ such 
that $\mathcal{V}=f^*V$.  This is immediate from the definition and can 
also be seen from the fibration 
$Gl(2n)/SO(2n) \to BSO(2n) \to BGl(2n)$ and the universal properties of 
$BSO(2n)$ and $BGl(2n)$.  

For any quadratic bundle $V \to X$, Edidin and Graham \cite{E&G} define a 
characteristic class $z_n$ in $CH^*(X)$ which maps to $2^{n-1}$ times the 
Euler class under the class map to cohomology.  For $X=Gl(2n)/SO(2n)$, we will 
compare this class with the class $y_n$ from Theorem 5.

We will begin by recalling their construction.  Let $Fl_i(V)$ be the 
variety of 
length $i$ isotropic flags in $V$.  A geometric point of $Fl_i(V)$ is a pair 
$(x, E_1 \subset E_2 \subset \cdots E_{i})$ where $x \in X$, and 
$E_j \subset V_x$ is an isotropic subspace of dimension $j$.  Write $Fl(V)$ 
for $Fl_{n-1}(V)$ and let 
$p:Fl(V) \to X$ and $p_i:Fl_i(V) \to Fl_{i-1}(V)$ be the obvious projection 
maps.  

We define a rank $2i$ vector bundle $V_{n-i}$ 
over $Fl_i(V)$ inductively.  Let $V_n=V$, the vector bundle over length 
zero isotropic flags.  $V_{n-1}$ will be the vector bundle whose fiber over 
$(x,E_1) \in Fl_1(V)$ is $E_1^{\perp}/E_1$.  This is well defined since 
$E_1$ is isotropic (i.e. $E_1\subset E_1^{\perp}$).  Given 
$V_{n-i-1}\to Fl_{i-1}(V)$, define $V_{n-i}$ to be the vector bundle whose 
fiber 
over $(x,E_1\subset\cdots\subset E_i)$ is $E_i^{\perp}/(E_i/E_{i-1})$ where 
$E_{i}^{\perp}$ is taken as a subvector space of 
${p_i^*V_{n-i-1}}_{(x,E_1\subset\cdots\subset E_i)}$.  Note that the 
$Fl_i(V)$ sit 
naturally in the projective space $\bP(V_{n-i-1})$ as the quadric of 
isotropic lines.  Let $h_{n-i}$ denote both the class of a hyperplane in 
$CH^1(\bP(V_{n-i-1}))$ and its pullback to $Fl_i(V)$.

Given 
$(x,E_1\subset\cdots\subset E_{n-1})\in Fl(V)$, $V_x$ has two rank $n$ 
isotropic 
sub-vectorspaces.  A choice of one of these determines a rank $n$ 
vectorbundle $V_n'$ since $det(V)$ was assumed to be trivial; this is a 
maximal isotropic subbundle of $p^*V$.  
The top Chern class of one 
of these bundles is the classical Euler class $c_n(V_n') \in H^*(Fl(V))$.  
Choice of $V_n'$ changes this by sign as $c_n(V')+c_n(V'')=0$ in the 
Chow ring for $X=Gl(2n)/SO(2n)$.
Define
$$(*) \hspace{1cm}z_n = p_*(h_2^2 \cdot h_3^4 \cdot ...\cdot h_n^{2n-2} \cdot c_n(V_n')).$$
to be the Chow Euler class in $CH^*(X)$.  This is equal to 
$$=p_*(h_2^3\cdot h_3^5\cdot ...\cdot h_n^{2n-1}\cdot c_1(V_1'))$$
where $V_1'$ and $V_1''$ are the two isotropic line subbundles of $V_1$.

\begin{theorem}The image of
Edidin-Graham's
element $z_n \in CH^n(BSO(2n))$ in $CH^n(GL(2n)/SO(2n))$ is the same as the 
non-zero class $y_n$
in $CH^n(Gl(2n)/SO(2n))$ from Corollary 3. 
\end{theorem}

\begin{proof}To compare these two classes, we will need to recall the 
definition of the class $y_n$ along with some of the proof of Theorem 5.  

Let $B$ be a Borel subgroup in $Gl(2n)$.  The $B$-orbits in $Gl(2n)/O(2n)$ 
are paramatrized by symmetric permutations in $2n$ letters, and those of 
$Gl(2n)/SO(2n)$ are components of the inverse images of $B$-orbits of the 
base. 

To calculate $CH^*(Gl(2n)/SO(2n))$, we decompose the space into $2n$ disjoint 
subvarieties, $X_i^{(n)}$, 
each of which is fibered over $Gl(2n-2)/SO(2n-2)$.  Define 
$(q,\epsilon)\in X^{(n)}_i$ if $q(e_{i+1},e_{2n})=\cdots=q(e_{2n},e_{2n})=0$ and 
$q(e_i,e_{2n})\neq 0$.  This is fibered over $Gl(2n-2)/SO(2n-2)$ via 
orthogonal projection onto $<e_{i},e_{2n}>^{\perp}$ where 
$<e_{i},e_{2n}>$ is the plane generated by the vectore $e_{i}$ and 
$e_{2n}$ on which the quadratic form is non-degenerate.  Call the projection onto the base $f^{(n)}_i$; this is a trivial 
fibration and the fiber is an 
open subset of affine space.  

The proof is by 
induction.  The Chow ring is generated by classes corresponding to the closure 
of $B$-orbits \cite{F-MacP-etc}, and there are only $3$ orbits in 
$Gl(2)/SO(2)$.  These are the dense orbit corresponding to the codimension 
zero cycle, and two closed orbits corresponding to codimension one 
subvarieties.  These are the two disjoint components of $X^{(1)}_1$, denoted 
$X^{(1)+}_1$ and $X^{(1)-}_1$.  
They are the inverse image from $Gl(2)/O(2)$ of a 
connected codimension one 
subvariety.  $Gl(2)/O(2)$ is isomorphic to an open subset of affine space, so 
by Lemma 1, it has trivial Chow ring.  Therefore 
$X^{(1)+}_1 + X^{(1)-}_1 \sim 0$, and $y_1$ is one of these two subvarieties.  

Assume $CH^*(Gl(2n-2)/SO(2n-2))=\Z \oplus \Z y_{n-1}$ and $y_{n-1}$ is the 
image of $X^{(1)+}_2$ via 
${f^{(n-1)}_{2n-3}}^*{f^{(n-2)}_{2n-5}}^*\cdots {f^{(2)}_3}^*$.  Since the 
$f^{(i)}_j$ are trivial fibrations with fibers an open subset of affine space, 
there is an induced surjection on Chow rings (\cite{Fulton} example 1.9.2).
Therefore $CH^*X^{(n)}_i$ is a quotient of $\Z\oplus \Z$.  In fact, two 
further 
induction arguments show that
$CH^*Gl(2n)/SO(2n)=\Z \oplus CH^{n-1}(X^{(n)}_{2n-1})=
\Z \oplus {f^{(n)}_{2n-1}}^*(y_{n-1})$.

We proceed with the proof of the theorem by induction on $n$.  For $n=0$ there is nothing to show.

Assume the case $n-1$.  Denote by
$\widetilde{V_{n-1}}\to Gl(2n-2)/SO(2n-2)$ the canonical bundle with rotating 
quadratic form. 
We assume that $z_{n-1}(\widetilde{V_{n-1}})$ is the class of the 
largest fixed point free orbit in $Gl(2n-2)/SO(2n-2)$, denoted $\widetilde{X}$.

In order to show $y_n(V_n)=p_*(h_2^3\cdots h_n^{2n-1}c_1(V_1'))=z_n(V_n)$, 
We will first show that $z_n$ is supported on $\overline{X_{2n-1}}$ ($y_n$ is supported there by definition).  
To do this, we will compute $h_n^{2n-1}$, where $h_n$ is the pullback of 
the hyperplane 
class from $\bP(V)$ to $Fl_1(V)$.  From there we will look at the 
restriction of $p_*(h_n^{2n-1})$ to $\widetilde{X}$.  

Since $V$ is a trivial vector bundle, 
$\bP(V)\cong X \times \bP(\C^{2n})$ as vector bundles, 
and the subvariety $h_n^{2n-1}$ in $\bP(V)$ is 
that of any point, take $(0:0:0:...:0:1)=e_{2n}$.
Since, 
$p_*(h_n^{2n-1})$ is the set of $q \in X$ such that the line generated by 
$e_{2n}$ is isotropic, it is the   
cycle representing ${\overline{X_{2n-1}}}$.  Therefore, $z_n$ is supported on this subvariety, and moreover $j^*(z_n)=f_{2n-1}^*(z_{n-1})$ where $j:X_{2n-1} \to \overline{X_{2n-1}}$ and $f_{2n-1}:X_{2n-1}\to Gl(2n-2)/SO(2n-2)$ is as previously defined.

As $j^*(y_n)=f_{2n-1}^*(y_{n-1})$, we have that 
$j^*(y_n-z_n)=0$.

By the fundamental exact sequence of Chow rings, this implies that $y_n-z_n$ is the pushforward 
of a cycle supported on $\overline{X_{2n-2}}$ as $\overline{X_{2n-2}}=\overline{X_{2n-1}} - X_{2n-1}$.

Since the fundamental exact sequence preserves dimension, $z_n-y_n$ must be a cycle of codimension $n-1$ in $\overline{X_{2n-2}}$ which itself has codimension $2$ in $Gl(2n)/SO(2n)$.  However, all cycles of codimension $n+1$ in $Gl(2n)/SO(2n)$ are rationally equivalent to zero, and therefore, $z_n=y_n$ in $CH^*(Gl(2n)/SO(2n))$.

\end{proof}

\section{$CH^*BSO(2n)$, Generators and Relations}
The main result of this paper is:
\begin{theorem}
$$CH^*BSO(2n) \cong 
\Z[c_2,c_3,...,c_2n,y_n]/(2c_{odd}=0, y_n^2=2^{2n-2}c_{2n}, y_n.c_{odd}=0).$$
\end{theorem}

\begin{proof}

By Theorems 4 and 5 and the fact that 
$$CH^*BGL(2n)=\Z[c_1,c_2,...,c_{2n}],$$
$CH^*BSO(2n)$ is a quotient of $\Z[c_1,...,c_{2n},y_n]$.  We begin by determining the relations.

The relation $c_1=0$ is simply the fact that the determinant line bundle 
$\Lambda^{2n}(V)$ for the 
standard representation of $SO(2n)$ is trivial.  

The relations $2c_{odd}=0$ hold because the standard 
representation of $SO(2n)$ is self-dual.

We now prove $y_n^2=(-1)^n2^{2n-2}c_{2n}$. 
Let $s=h_2^2 h_3^4\cdots h_n^{2n-1}$.
We have 
$p_*(s.p^*x)=2^{n-1}x$ for any $x \in CH^*X$\cite{E&G}.  
Now
\\
$p_*(s.p^*(c_{2n}(V)))=2^{n-1}c_{2n}(V)$, so\\
$p_*(s.p^*(p_*(s.p^*(c_{2n}(V)))))=2^{2n-2}c_{2n}(V) \mbox{ in } CH^*X$.  \\

Since
$p^*(c_{2n}(V))=c_{2n}(p^*V)=c_n(V')^2 \mbox{ in } CH^* Fl(V)$,
where $V'$ is a maximal isotropic subbundle of $p^*V$, using the projection formula twice:\\
$$
\begin{array}{cl}
p_*(s.p^*(p_*(s.c_n(V')^2)))&=2^{2n-2}c_{2n}(V)\\
\quad&=2^{n-1}.p_*(s.c_n(V')^2)\\
\quad&=p_*(s.c_n(V').2^{n-1}c_n(V'))\\
\quad&=p_*(s.c_n(V').p^*(y_n))\\
\quad&=y_n^2.
\end{array}
$$

Finally, we must prove $y_n.c_{odd}=0$. 
It is enough to show 
$c_n(V').c_{odd}(p^*V)=0$ in $CH^*(Fl(V)$, as
$p_*(s.c_n(V').c_{odd}(p^*V))=y_nc_{odd}(V)$.
Over $Fl(V)$ 
we have an exact sequence of vector bundles, where $V'$ is a maximal 
isotropic subbundle of $p^*(V)$ and $V'^\vee$ is its dual.
$$0 \to V' \to p^*V \to V'^\vee \to 0.$$
It follows that
$c(p^*V)=c(V').c(V'^\vee)$
and hence the odd Chern classes of $p^*V$ are 
zero.

This shows that 
$CH^*BSO(2n)$ is a quotient of the ring 
$$\Z [c_2,...,c_{2n},y_n]/(2c_{odd}=0, y_n^2=2^{2n-2}c_{2n}, y_n.c_{odd}=0).$$

The computation of $H^*BSO(2n)$, see for example \cite{M&S} gives 
$$H^*BSO(2n) \supset \Z[c_2,c_3,..,c_{2n},e]/(2c_{odd}=0, e^2=c_{2n})$$   
\noindent where $e$ is the Euler class of the 
standard representation.  We have shown that the composite 
$$
\Z[c_2,...,c_{2n},y_n]/(2c_{odd},y_n.c_{odd},y_n^2-2^{2n-2}c_{2n})
\to CH^*BSO(2n) \to H^*BSO(2n)
$$
is injective and that the first map is surjective.  Therefore, the first map 
is an isomorphism.
\end{proof}

As an immediate corollary, we have 

\begin{corollary}
$CH^*(BSO(2n)) \hookrightarrow H^*(BSO(2n)).$
\end{corollary}

The relation $y_n.c_{odd}=0$ deserves a short comment.
Unlike the other relations, the corresponding relation 
in cohomology is trivial, as  
$y_n.c_{odd}$ maps to $2^{n-1}\chi.c_{odd}$, and 
odd Chern classes are $2$-torsion (for the same reason they are $2$-torsion in the Chow ring). 
In $CH^*BSO(2n)$, however, 
the class $y_n$ is not $2.\mbox{ anything}$.  (Even in $CH^*BSO(4,\C)$ this 
is apparent, according to Pandharipande's computation \cite{Rahul} 
$y_2=f_2-c_2$ where $f_2$ 
is the second Chern class of the exceptional representation of $SO(4)$.)  
So unlike the other relations, looking to cohomology is not productive.

\section{Stratifying quadratic bundles}

The space of quadratic bundles is naturally stratified.  A quadratic bundle 
is locally trivial in the Etale topology \cite{Swan}, but not in the Zariski 
topology.  In this section we stratify the space of quadratic 
bundles by ``how Zariski locally trivial they are''.  The closed strata is 
the set of ZLT quadratic bundles.  The generic strata (which contains the 
universal trivial bundle over $Gl(2n)/SO(2n)$) is furthest from this.  
We will see that the strata a bundle is in determines the two divisibility of its Chow Euler class in cohomology.  We begin with a definition.

\begin{definition}
We say that a rank $N$ vector bundle $V \to X$ with non-degenerate quadratic 
form is {\bf Zariski locally trivial to codimension k} if there 
exists a Zariski open cover $\{U_i\}$ of $X$, such that 
$V\mid_{U_i}\simeq X \times A^N$ and, under these isomorphisms, the map 
$U_i \to Gl(N)/SO(N)$ has image in the subgroup
$$
\left(\begin{array}{cccc}
Gl(\mbox{k})/SO(\mbox{k})&\quad&\quad&\quad\\
\quad&\quad&\quad&1\\
\quad&\quad&\cdots&\quad\\
\quad&1&\quad&\quad
\end{array}\right) 
\mbox{or}
\left(\begin{array}{cccc}
Gl(\mbox{k})/SO(\mbox{k})&\quad&\quad&\quad\\
\quad&\quad&\quad&1\\
\quad&\quad&\cdots&\quad\\
\quad&1&\quad&\quad\\
\quad&\quad&\quad&f(u)
\end{array}\right)
$$
where 
$k$ is even and odd 
respectively and $f(u)$ is 
a nowhere vanishing function on $U_i$.  
\end{definition}

In other words, a quadratic
bundle $V \to X$ is Zariski locally trivial to codim $2k$ if there 
exists a Zariski trivialization mapping to the universal trivial quadratic bundle as follows:
$$
\xymatrix{
V\mid_{U_i}\ar[d]\ar[r] & Gl(2n)/SO(2n)\times C^{2n}\ar[d]\\
X\mid_{U_i}\ar[r] & Gl(2n)/SO(2n)\\
}
$$
whose image in $Gl(2n)/SO(2n)$ is of the above form.

Notice we are not requiring any compatibility on the intersections $U_i \cap U_j$.

\begin{proposition}
If $V\to X$ is a 
quadratic bundle
which is Zariski locally trivial to codimension $2k$ for $1 \leq k \leq n$, 
then there is 
a class $x$ in the Chow ring of $X$ mapping to $2^{k-1}\chi$ in cohomology. 
\end{proposition}
\begin{proof}

Recall a proper map $f:X' \to X$ is a Chow envelope in the sense of \cite{Fulton} if for every subvariety $Z\subset X$, there is a subvariety $Z' \subset X'$ mapping birationally onto $Z$.  In that case, $f_*:CH^*X' \to CH^*X$ is surjective, $f^*:CH^*X\to CH^*X'$ is injective, and $\alpha \in CH^*X'$ is the pullback of a class in $CH^*X$ if and only if it satisfies the descent condition $\pi_1^*\alpha=\pi_2^*\alpha$, where $\pi_i:X'\times_X X'\to X'$ are the two projection maps\cite{Kimura}.  Moreover, the property of being a Chow envelope is clearly Zariski local in $X$.

Now apply this to $p_1:Fl_{n-k}(V)\to X$.  To see this is a Chow envelope, it is enough to check on a chart in a cover as in Definition 1, i.e. when $V=X\times \C^{2n}$ and the quadratic form on $\C^{2n}$ has a length $n-k$ isotropic flag $F_{n-k}=\langle e_{2n}\rangle \subset\langle e_{2n},e_{2n-1}\rangle\subset \cdots\subset \langle e_{2n},...,e_{n+k+1}\rangle$.  Given a subvariety $Z\subset X$, the subvariety $Z'=\{(z,F_{n-k}\mid z\in Z\}\subset Fl_{n-k}(V)$ maps birationally onto $Z$.

Set $\tilde{x}={p_2}_*(h_2^2\cdots h_k^{2k-2}\cdot c_n(V'))\in Fl_{n-k}(V)$, where $p_2:Fl_{n-1}(V)\to Fl_{n-k}(V)$.  It is clear that $\tilde{x}$ satisfies the descent condition, and moreover maps to $2^{k-1}\chi$ in $H^*(Fl_{n-k}(V))$.  Therefore, there exists an element $x \in CH^*X$ such that $p_1^*(x)=\tilde{x}$ and $\sc{cl}(x)=2^{k-1}\chi\in H^*X$.

\end{proof}

Note that every quadratic bundle which is Zariski locally trivial to codimension $2$
is, in fact, Zariski locally trivial.

\section{Some Representation Theory}
Now we ask whether our class $y_n$ mapping to $2^{n-1}\chi$ in cohomology 
is a polynomial in  
Chern classes of representations of $SO(2n)$.  The representation ring 
of $SO(2n)$ 
(homotopy equivalent to its compact form) 
is generated by 
the standard $2n$ dimensional representation $V$, exterior powers of $V$, 
and one additional irreducible 
representation. 
$D_n^+$ is the 
irreducible representation of $SO(2n)$ whose highest weight vector is twice 
that of either of the two half spin representations of $\mathfrak{so}_{2n}$; 
it can be realized as the 
$+1$ eigenspace of the Hodge star operator 
on $\Lambda^nV$.
If $y_n$ is to be a polynomial in 
Chern classes of representations, then it must be  
$y_n = \pm c_n(D_n^+)$ mod Chern classes of the standard representation.  
However,
we have the following result.

\begin{theorem} Let $D^+_n$ be the space of self-dual n-forms, as above.  
Then in $H^*BSO(2n)$, 
$c_n(D^+_n) = \pm (n-1)!2^{n-1}y_n + p(c_1,c_2,...,c_{2n})$ 
\noindent 
where $p(c_1,c_2,...c_{2n})$ is a polynomial in the Chern classes of the 
standard representation of $SO(2n)$.
\end{theorem}

This gives the immediate corollary:

\begin{corollary}
The ring $CH^*BSO(2n)$ is not generated by Chern classes of any 
representation for $n \geq 3$.
\end{corollary}

As Totaro pointed out, this corollary gives the 
first known  
example
of a
compact Lie group $G$ and a prime number $p$ such that $H^*(BG,Z/p)$ is
concentrated in even degrees (known since the 1950's), but is not 
generated by Chern classes of
representations of $G$.  
For example, take $G=SO(8)$ and $p=3$.
Then $CH^*(BSO(8))/3$ maps isomorphically to
$H^*BSO(8,Z/3)$.  Similar properties hold for $BSO(2n)$ for $n\geq 5$ and 
$p$ any prime $2<p\leq n$.

\begin{proof}

Since 
$$
CH^*(BSO(2n))\otimes \Q \hookrightarrow H^*(BSO(2n),\Q) \cong H^*(BSO(2n,\R),\Q) \hookrightarrow H^*(B(S^1)^n,\Q),
$$
we can compute the coefficient of $y_n$ in $c_n(D_n^+)$ by restricting to the maximal torus.
Let $z_i$ be the representation of 
$(S^1)^n$ that restricts to the $i^{th}$ factor. 

From the definition of $D_n^+$ as the space of self dual $n$-forms, the total Chern class
$$
ch(D_n^+)=f(z_1, z_2, \dots, z_n) = \prod (1 + \sum e_i z_i)
$$
where the product is taken over all choices of $\{e_i=\pm 1\}$ with
an even number of $+1$'s.
We will compute the coefficient of the Euler class $z_1 z_2 \dots z_n$ in this polynomial.

The following combinatorial argument was communicated to me by Sandy Kutin.

Note that $f(z_1, z_2, \dots, z_n)$ = $f(a_1 z_1, a_2 z_2, \dots, a_n z_n)$
for any sequence of $a_i=\pm 1$ such that $\prod a_i = 1$, as  
this simply permutes the terms in the product.
In particular, the contribution to the monomial $z_1 z_2 \dots z_n$
is the same regardless of which term the $z_n$ comes from.
So we can compute by choosing a particular term for the $z_n$ to come from 
and multiply by the total number of terms.  We choose the term 
$(1-z_1-\cdots-z_n)$ and note that there are $2^{n-1}$ terms.

Since we have accounted for the $z_n$ term, this contribution can be
written as 
(the negative of) 
the coefficient of $z_1 z_2 \dots z_{n-1}$ in
$$
g(z_1, z_2, \dots, z_{n-1}) = \prod (1 + \sum e_i z_i)
$$
where the product is now taken over all choices of $\{e_i = \pm 1\}$
except for $1 - \sum z_i$ since any
one of these can be extended to a term in the original product by choosing 
the coefficient $e_n=\pm 1$ appropriately.

Note that the product of the terms $(1+e_1z_1+\cdots + e_{n-1}z_{n-1})$ and 
$(1+e_1z_1+\cdots +e_{i-1}z_{i-1}-e_iz_i+e_{i+1}z_{i+1}+\cdots +e_{n-1}z_{n-1})$ 
has no monomials whose power of $z_i$ is less than $2$.  Therefore, such a 
pair can not contribute a $z_i$ to the monomial $z_1z_2...z_{n-1}$.  The 
proof centers around detecting such pairs.  

We first look for which terms can contribute a $z_1$.  Since the product 
$g(z_1,...,z_{n-1})$ was taken over all choices of $\{e_i=\pm1\}$ except 
$1-\sum e_iz_i$, the only term which will not have a pair of the above type is 
$1+z_1-z_2-\cdots -z_{n-1}$.  Therefore, the $z_1$ in the monomial 
$z_1z_2\cdots z_{n-1}$ must come from this term.

Now look for which terms can contribute the $z_2$.  There will be two terms 
that do not pair in the above way.  These are $1+z_1+z_2-z_3-\cdots -z_{n-1}$ 
whose pair has already contributed its $z_1$ term, and 
$1-z_1+z_2-z_3-\cdots -z_{n-1}$ whose pair had already been removed from $g$ 
after contributing its $z_n$ term.  We will pick the first one and multiply 
by $2$.  

In general, when looking for terms that can contribute their $z_i$ term to 
the monomial $z_1z_2\cdots z_{n-1}$, we will find there are $i$ factors that 
remain unpaired.  These are $(1+z_1+\cdots+z_i-z_{i+1}-\cdots -z_{n-1}), 
(1+z_1+\cdots +z_{i-2}-z_{i-1}+z_i-z_{i+1}-\cdots -z_{n-1}),\cdots ,
(1+z_1-z_2-\cdots-z_{i-1}+z_i-z_{i+1}-\cdots -z_{n-1}),
(1-z_1-\cdots -z_{i-1}+z_i-z_{i+1}-\cdots -z_{n-1}.$  Again, pick the top 
term and multiply by $i$.  

Since we go through this process $n-1$ times, in total, we will get a 
coefficient of $(n-1)!$ for the $z_1\cdots z_{n-1}$ term of the 
polynomial $g$, so a coefficient of 
$2^{n-1}(n-1)!$ for the monomial $z_1\cdots z_n$ in the 
original polynomial $f$.

\end{proof}

\bibliographystyle{plain}

\end{document}